\documentclass[12pt,reqno]{amsart}
\usepackage{enumerate, latexsym, amsmath, amsfonts, amssymb, amsthm, color}
\allowdisplaybreaks[4]
\def\pmod #1{\ ({\rm{mod}}\ #1)}

\def\bg{\bigg}
\def\({\bg(}
\def\){\bg)}

\def\Ack{\medskip\noindent {\bf Acknowledgment}}
\theoremstyle{plain}
\newtheorem{theorem}{Theorem}

\newtheorem{lemma}{Lemma}

\newtheorem{conjecture}{Conjecture}
\theoremstyle{definition}

\theoremstyle{remark}

 \vspace{4mm}

\begin{document}
\title
[{Proof of two congruence conjectures of Z.-W. Sun}]
{Proof of two congruence conjectures of Z.-W. Sun}

\author
[Guo-Shuai Mao] {Guo-Shuai Mao}

\address {(Guo-Shuai Mao) Department of Mathematics, Nanjing
University of Information Science and Technology, Nanjing 210044,  People's Republic of China\\
{\tt maogsmath@163.com  } }

\keywords{Congruences; harmonic numbers; Lucas sequence; Legendre symbol; binary quadratic forms.
\newline \indent This work was supported by the National Natural Science Foundation of China (grant 12001288).}

 \begin{abstract} In this paper, we mainly prove two congruence conjecture of Z.-W. Sun. Let $p\equiv3\pmod 4$ be a prime. Then
$$\sum_{k=0}^{p-1}\frac{\binom{2k}k^2}{8^k}\equiv-\sum_{k=0}^{p-1}\frac{\binom{2k}k^2}{(-16)^k}\pmod{p^3}.$$
And for any odd prime $p$, if $p=x^2+y^2$ with $4|x-1, 2|y$, then
$$
\sum_{k=0}^{p-1}\frac{(k+1)\binom{2k}k^2}{8^k}+\sum_{k=0}^{(p-1)/2}\frac{(2k+1)\binom{2k}k^2}{(-16)^k}\equiv2\left(\frac{2}p\right)x\pmod{p^3}.
$$
\end{abstract}

\maketitle

\section{Introduction}
\setcounter{lemma}{0}
\setcounter{theorem}{0}
\setcounter{corollary}{0}
\setcounter{remark}{0}
\setcounter{equation}{0}
\setcounter{conjecture}{0}
Let $p>3$ be a prime. Rodriguez-Villegas \cite{RV} conjectured that
$$
\sum_{k=0}^{(p-1)/2}\frac{\binom{2k}k^2}{16^k}\equiv\left(\frac{-1}{p}\right)\pmod{p^2},
$$
which was confirmed by E. Mortenson \cite{Mo1} via an advanced tool involving the $p$-adic Gamma function and the Gross-Koblitz formula for character sums, the above $\left(\frac{\cdot}p\right)$ stands for the Legendre Symbol. Sun \cite{sffa} generalized the above congruence to the following type: Let $p>3$ be a prime and let $d\in\{0,1,\ldots,(p-1)/2\}$. Then
$$
\sum_{k=0}^{(p-1)/2}\frac{\binom{2k}k\binom{2k}{k+d}}{16^k}\equiv\left(\frac{-1}{p}\right)+p^2\frac{(-1)^d}{4}E_{p-3}\left(d+\frac12\right)\pmod{p^3}.
$$
This, with $d=0$ yields that
$$
\sum_{k=0}^{(p-1)/2}\frac{\binom{2k}k^2}{16^k}\equiv\left(\frac{-1}{p}\right)+p^2E_{p-3}\pmod{p^3},
$$
which was first proved by Sun \cite{sscm} with the help of the software \texttt{Sigma} \cite{S}. For $\mathbb{N}=\{0,1,2,\ldots\}$, the above $\{E_n\}_{n\in\mathbb{N}}$ and $\{E_n(x)\}_{n\in\mathbb{N}}$ are Euler numbers and Euler polynomials defined by
\begin{gather*}
E_0=1\ \ \ \mbox{and}\ \ \sum_{\substack{k=0\\2|k}}^n\binom{n}kE_{n-k}=0\ \ \ \mbox{for}\ \  n\in\mathbb{Z}^{+}=\{1,2,3,\ldots\},\\
E_n(x)=\sum_{k=0}^n\binom{n}k\frac{E_k}{2^k}\left(x-\frac12\right)^{n-k}.
\end{gather*}
In \cite{sffa}, Sun also showed that if $p\equiv3\pmod4$ is a prime, then
$$
\sum_{k=0}^{(p-1)/2}\frac{\binom{2k}k^2}{8^k}\equiv-\sum_{k=0}^{(p-1)/2}\frac{\binom{2k}k^2}{(-16)^k}\equiv\frac{(-1)^{(p+1)/4}2p}{\binom{(p+1)/2}{(p+1)/4}}\pmod{p^2}.
$$
In \cite[Conjecture 5.5]{sjnt}, Sun proposed the following conjecture.
\begin{conjecture}\label{Conj1.1} If $p\equiv3\pmod4$ is a prime, then
\begin{equation}\label{zhu}
 \sum_{k=0}^{p-1}\frac{\binom{2k}k^2}{8^k}\equiv-\sum_{k=0}^{p-1}\frac{\binom{2k}k^2}{(-16)^k}\pmod{p^3}.
 \end{equation}
\end{conjecture}
\noindent Our first goal in this paper is to prove the above conjecture.
 \begin{theorem}\label{Th1.1} Conjecture \ref{Conj1.1} is true.
\end{theorem}
\noindent Let $p$ be an odd prime. Sun \cite[Remark 1.1]{saa} conjectured that
\begin{align*}
&\sum_{k=0}^{p-1}\frac{(k+1)\binom{2k}k^2}{8^k}+\sum_{k=0}^{(p-1)/2}\frac{(2k+1)\binom{2k}k^2}{(-16)^k}\\
&\equiv\begin{cases} 2\left(\frac{2}p\right)x\ \pmod{p^3} &\mbox{if}\ \it{p}={\textit{x}}^{\texttt{2}}+\textit{y}^{\texttt{2}}\ (\texttt{4}|(\textit{x}-\texttt{1})\ \& \ \texttt{2}|\textit{y}), \\
0\ \pmod{p^2}  &\mbox{if}\ \it{p}\equiv\texttt{3}\pmod {\text{4}}.\end{cases}
\end{align*}
The case $p\equiv3\pmod4$ of this conjecture was confirmed by Sun in \cite[Theorem 1.3(i)]{sffa}. Now we prove the case $p\equiv1\pmod4$ of this conjecture.
\begin{theorem}\label{Th1.2} For any odd prime $p$, if $p=x^2+y^2$ with $4|x-1, 2|y$, then
$$
\sum_{k=0}^{p-1}\frac{(k+1)\binom{2k}k^2}{8^k}+\sum_{k=0}^{(p-1)/2}\frac{(2k+1)\binom{2k}k^2}{(-16)^k}\equiv2\left(\frac{2}p\right)x\pmod{p^3}.
$$
\end{theorem}
The paper proceeds as follows. We are going to prove Theorem \ref{Th1.1} in the next section. Section 2 is devoted to proving Theorem \ref{Th1.2}. Our proofs make use of some congruences involving harmonic numbers, combinatorial identities, in parts supported by the symbolic summation package \texttt{Sigma} \cite{S}.
\section{Proof of Theorem \ref{Th1.1}}
 \setcounter{lemma}{0}
\setcounter{theorem}{0}
\setcounter{corollary}{0}
\setcounter{remark}{0}
\setcounter{equation}{0}
\setcounter{conjecture}{0}
For $n,m\in\mathbb{Z}^{+}$, define
$$
H_n^{(m)}:=\sum_{1\leq k\leq n}\frac1{k^m},\ \ H_0^{(m)}:=0,
$$
these numbers with $m=1$ are often called the classic harmonic numbers.
\begin{lemma}\label{Lem2.0}{\rm(\cite{s2000,s2008})} Let $p>5$ be a prime. Then
\begin{gather*}
H_{\lfloor p/4\rfloor}\equiv-3q_p(2)+\frac32pq^2_p(2)-(-1)^{(p-1)/2}pE_{p-3}\pmod{p^2},\\
H_{(p-1)/2}\equiv-2q_p(2)\pmod{p},\ \ H_{p-1}\equiv0\pmod{p^2},\\
H_{p-1}^{(2)}\equiv H_{(p-1)/2}^{(2)}\equiv0\pmod{p},\ \ H_{\lfloor p/4\rfloor}^{(2)}\equiv(-1)^{(p-1)/2}4E_{p-3}\pmod p,
\end{gather*}
where $q_p(a)=(a^{p-1}-1)/p$ denotes the Fermat quotient.
\end{lemma}
By the above lemma, it is easy to see that for any odd prime $p$ and each integer
$0\leq k\leq p-1$,
\begin{equation}\label{kp-1-k}
H_{p-1-k}=\sum_{j=1}^{p-1-k}\frac1j=\sum_{j=k+1}^{p-1}\frac1{p-j}\equiv H_k-H_{p-1}\equiv H_k\pmod p.
\end{equation}
\begin{lemma}\label{Lem2.1}{\rm (\cite[Lemma 4.2]{sijm})} Let $p=2n+1$ be an odd prime, and let $k\in\{0,1,\ldots,n\}$. Then
\begin{align*}
\frac{\binom{n}k}{\binom{2k}k/(-4)^k}\equiv&1-p\sum_{j=1}^k\frac1{2j-1}+\frac{p^2}2\left(\sum_{j=1}^k\frac1{2j-1}\right)^2\\
&-\frac{p^2}2\sum_{j=1}^k\frac1{(2j-1)^2}\pmod{p^3}.
\end{align*}
\end{lemma}
\begin{lemma}\label{mpt} Let $p>5$ be a prime with $p\equiv3\pmod4$. Then
\begin{equation*}
 \binom{-1/4}{\frac{p-1}2}^{-1}\equiv\frac{(-1)^{(p+1)/4}}{\binom{\frac{p-1}2}{\frac{p-3}4}}(1-3p+pq_p(2))\pmod{p^2}.
 \end{equation*}
 \end{lemma}
\begin{proof} Set $m=(3p-1)/4, t=-3/4$. By Lemma \ref{Lem2.0} and (\ref{kp-1-k}), we have
\begin{align*}
\binom{m+pt}{(p-1)/2}&=\frac{(m+pt)\cdots(m+pt-(p-1)/2+1)}{((p-1)/2)!}\\
&\equiv\frac{m\cdots(m-(p-1)/2+1)}{((p-1)/2)!}(1+pt(H_m-H_{m-(p-1)/2}))\\
&\equiv\binom{(3p-1)/4}{(p-1)/2}(1+3p)\pmod{p^2}.
\end{align*}
And
\begin{align*}
\binom{(3p-1)/4}{(p-1)/2}&=\binom{(3p-1)/4}{(p+1)/4}=(-1)^{(p+1)/4}\binom{-p+(p-1)/2}{(p+1)/4}\\
&\equiv(-1)^{(p+1)/4}\binom{(p-1)/2}{(p+1)/4}(1-p(H_{(p-1)/2}-H_{(p-3)/4}))\\
&\equiv(-1)^{(p+1)/4}\binom{(p-1)/2}{(p+1)/4}(1-pq_p(2))\pmod{p^2}.
\end{align*}
We immediately obtain that
$$\binom{-1/4}{\frac{p-1}2}\equiv(-1)^{(p+1)/4}\binom{\frac{p-1}2}{\frac{p-3}4}(1+3p-pq_p(2))\pmod{p^2},$$
hence
$$\binom{-1/4}{\frac{p-1}2}^{-1}\equiv\frac{(-1)^{(p+1)/4}}{\binom{\frac{p-1}2}{\frac{p-3}4}}(1-3p+pq_p(2))\pmod{p^2}.$$
Now the proof of Lemma \ref{mpt} is completed.
\end{proof}
Let $p>3$ be a prime. In 1862, Wolstenholme \cite{Wol} established the well-known congruence
\begin{equation}\label{wol}
\frac12\binom{2p}p=\binom{2p-1}{p-1}\equiv1\pmod{p^3}.
\end{equation}
\noindent In view of \cite[(3.99)]{G}, we have the following identity,
\begin{equation}\label{id0}
\binom{2n}n=\sum_{k=0}^{\lfloor n/2\rfloor}\binom{2k}k\binom{n}{2k}2^{n-2k}.
\end{equation}
Now we evaluate $\sum_{k=0}^{p-1}\binom{2k}k^2/8^k$ modulo $p^3$. By (\ref{id0}), we have
$$
\sum_{k=0}^{p-1}\frac{\binom{2k}k^2}{8^k}=\sum_{k=0}^{p-1}\frac{\binom{2k}k}{8^k}\sum_{j=0}^{\lfloor k/2\rfloor}\binom{2j}j\binom{k}{2j}2^{k-2j}=\sum_{j=0}^{(p-1)/2}\frac{\binom{2j}j}{4^j}\sum_{k=2j}^{p-1}\frac{\binom{2k}k\binom{k}{2j}}{4^k}.
$$
By \texttt{Sigma}, we can find and prove the following identity:
$$
\sum_{k=2j}^{n-1}\frac{\binom{2k}k\binom{k}{2j}}{4^k}=\frac{n\binom{2n-1}{n-1}\binom{n-1}{2j}}{4^{n-1}(4j+1)}.
$$
Substituting $n=p$ into the above identity and by (\ref{wol}) we have
\begin{align*}
\sum_{k=0}^{p-1}\frac{\binom{2k}k^2}{8^k}&=\sum_{j=0}^{(p-1)/2}\frac{\binom{2j}j}{4^j}\frac{p\binom{2p-1}{p-1}\binom{p-1}{2j}}{4^{p-1}(4j+1)}\\
&\equiv \frac{p}{4^{p-1}}\sum_{j=0}^{(p-1)/2}\frac{\binom{2j}j}{4^j}\frac{1-pH_{2j}}{4j+1}\pmod{p^3},
\end{align*}
where we also used the Fermat's Little Theorem and the fact that $$\binom{p-1}{2j}=\prod_{i=1}^{2j}\left(\frac{p-i}i\right)=\prod_{i=1}^{2j}\left(1-\frac{p}i\right)\equiv1-pH_{2j}\pmod{p^2}.$$
In view of Lemma \ref{Lem2.1}, we have
\begin{align*}
\sum_{j=0}^{\frac{p-1}2}\frac{\binom{2j}j}{(4j+1)4^j}\equiv\sum_{j=0}^{\frac{p-1}2}\frac{\binom{\frac{p-1}2}j(-1)^j}{(4j+1)}\left(1+p\left(H_{2j}-\frac12H_j\right)\right)\pmod{p^2}.
\end{align*}
Hence
$$
\sum_{k=0}^{p-1}\frac{\binom{2k}k^2}{8^k}\equiv \frac{p}{4^{p-1}}\sum_{j=0}^{\frac{p-1}2}\frac{\binom{\frac{p-1}2}j(-1)^j}{(4j+1)}-\frac{p^2}2\sum_{j=0}^{\frac{p-1}2}\frac{\binom{\frac{p-1}2}j(-1)^j H_j}{(4j+1)}\pmod{p^3}.
$$
By \texttt{Sigma} again, we can find and prove the following identities:
\begin{gather}
\sum_{k=0}^n\frac{\binom{n}k(-1)^k}{4k+1}=\frac1{4n+1}\frac{(-1)^n}{\binom{-1/4}n},\label{nk4k+1}\\
\sum_{k=0}^n\frac{\binom{n}k(-1)^k}{4k+1}H_k=-\frac1{4n+1}\frac{(-1)^n}{\binom{-1/4}n}\sum_{k=1}^n\frac{\binom{-1/4}k(-1)^k}{k}.\label{nk4k+1hk}
\end{gather}
Setting $n=(p-1)/2$ into the above identities and noting that $n\equiv1\pmod2$, we have
$$
\sum_{k=0}^{p-1}\frac{\binom{2k}k^2}{8^k}\equiv \frac{p}{4^{p-1}}\frac{1+2p}{\binom{-1/4}{(p-1)/2}}+\frac{p^2}2\frac{1}{\binom{-1/4}{(p-1)/2}}\sum_{k=1}^{(p-1)/2}\frac{\binom{-1/4}k(-1)^k}{k}\pmod{p^3}.
$$
In view of Lemma \ref{mpt} and note that $4^{1-p}=(1+pq_p(2))^{-2}\equiv1-2pq_p(2)\pmod{p^2}$, we have
$$
\sum_{k=0}^{p-1}\frac{\binom{2k}k^2}{8^k}\equiv \frac{(-1)^{\frac{p+1}4}}{\binom{\frac{p-1}2}{\frac{p-3}4}}\left(p(1-p-pq_p(2))+\frac{p^2}2\sum_{k=1}^{\frac{p-1}2}\frac{\binom{-\frac14}k(-1)^k}{k}\right)\pmod{p^3}.
$$
By \cite[(1.134)]{G} and note that $(p+1)/2\equiv0\pmod2$, we have
\begin{align*}
\sum_{k=\frac{p+1}2}^{\frac{3p-1}4}\frac{\binom{\frac{3p-1}4}k(-1)^k}{k}&=(-1)^{\frac{p+1}2}\sum_{k=\frac{p+1}2}^{\frac{3p-1}4}\frac{\binom{k-1}{\frac{p-1}2}}{k}=\sum_{k=0}^{ (p-3)/4}\frac{\binom{k+(p-1)/2}{(p-1)/2}}{k+(p+1)/2}\\
&\equiv2\sum_{k=0}^{(p-3)/4}\frac{\binom{k+(p-1)/2}{k}}{2k+1}\equiv2\sum_{k=0}^{(p-3)/4}\frac{\binom{-\frac{p+1}2}{k}(-1)^k}{2k+1}\\
&\equiv2\sum_{k=0}^{(p-3)/4}\frac{\binom{\frac{p-1}2}{k}(-1)^k}{2k+1}\pmod p.
\end{align*}
Substituting $n=(3p-1)/4$ into \cite[(1.45)]{G}, then by (\ref{kp-1-k}) and Lemma \ref{Lem2.0}, we have
\begin{align*}
\sum_{k=1}^{\frac{p-1}2}\frac{\binom{-\frac14}k(-1)^k}{k}&\equiv\sum_{k=1}^{\frac{3p-1}4}\frac{\binom{\frac{3p-1}4}k(-1)^k}{k}-\sum_{k=\frac{p+1}2}^{\frac{3p-1}4}\frac{\binom{\frac{3p-1}4}k(-1)^k}{k}\\
&=-H_{(3p-1)/4}-2\sum_{k=0}^{(p-3)/4}\frac{\binom{\frac{p-1}2}{k}(-1)^k}{2k+1}\\
&\equiv3q_p(2)-2\sum_{k=0}^{(p-3)/4}\frac{\binom{\frac{p-1}2}{k}(-1)^k}{2k+1}\pmod p.
\end{align*}
Therefore, modulo $p^3$, we have
\begin{align}\label{zhu1}
&\sum_{k=0}^{p-1}\frac{\binom{2k}k^2}{8^k}\notag\\
&\equiv \frac{(-1)^{\frac{p+1}4}}{\binom{\frac{p-1}2}{\frac{p-3}4}}\left(p\left(1-p+\frac12pq_p(2)\right)-p^2\sum_{k=0}^{(p-3)/4}\frac{\binom{\frac{p-1}2}{k}(-1)^k}{2k+1}\right).
\end{align}
\begin{lemma}\label{Lem2.2} For any prime $p>5$ and each integer $0\leq k\leq(p-1)/2$, we have
\begin{gather*}
H_{p-1-2k}\equiv H_{2k}+pH_{2k}^{(2)}\pmod{p^2},\\
H_{(p-1)/2-k}\equiv H_{(p-1)/2}+2H_{2k}-H_k+2pH_{2k}^{(2)}-\frac{p}2H_k^{(2)}\pmod{p^2},\\
H_{p-1-2k}^{(2)}\equiv-H_{2k}^{(2)}\pmod p,\ \ H_{(p-1)/2-k}^{(2)}\equiv H_k^{(2)}-4H_{2k}^{(2)}\pmod p.
\end{gather*}
\end{lemma}
\begin{proof} In view of Lemma \ref{Lem2.0}, we have
\begin{align*}
H_{p-1-2k}&=\sum_{j=1}^{p-1-2k}\frac1j=\sum_{j=2k+1}^{p-1}\frac1{p-j}=\sum_{j=2k+1}^{p-1}\frac{p+j}{p^2-j^2}\equiv-\sum_{j=2k+1}^{p-1}\frac{p+j}{j^2}\\
&=-(pH_{p-1}^{(2)}-pH_{2k}^{(2)}+H_{p-1}-H_{2k})\equiv H_{2k}+pH_{2k}^{(2)}\pmod{p^2},
\end{align*}
\begin{align*}
&H_{(p-1)/2-k}=\sum_{j=1}^{\frac{p-1}2-k}\frac1j=\sum_{j=k+1}^{\frac{p-1}2}\frac1{\frac{p+1}2-j}\\
&=2\sum_{j=k+1}^{\frac{p-1}2}\frac{p+2j-1}{p^2-(2j-1)^2}\equiv-2\sum_{j=k+1}^{\frac{p-1}2}\frac{p+2j-1}{(2j-1)^2}\\
&=-2\left(p\left(H_{p-1}^{(2)}-\frac14H_{\frac{p-1}2}^{(2)}-H_{2k}^{(2)}+\frac14H_k^{(2)}\right)+H_{p-1}-\frac12H_{\frac{p-1}2}-H_{2k}+\frac12H_k\right)\\
&\equiv H_{(p-1)/2}+2H_{2k}-H_k+2pH_{2k}^{(2)}-\frac{p}2H_k^{(2)}\pmod{p^2},
\end{align*}
\begin{align*}
H_{p-1-2k}^{(2)}&=\sum_{j=1}^{p-1-2k}\frac1{j^2}=\sum_{j=2k+1}^{p-1}\frac1{(p-j)^2}\equiv\sum_{j=2k+1}^{p-1}\frac{1}{j^2}\\
&=H_{p-1}^{(2)}-H_{2k}^{(2)}\equiv -H_{2k}^{(2)}\pmod{p}
\end{align*}
and
\begin{align*}
&H_{(p-1)/2-k}^{(2)}=\sum_{j=1}^{\frac{p-1}2-k}\frac1{j^2}=\sum_{j=k+1}^{\frac{p-1}2}\frac1{\frac{(p+1}2-j)^2}\equiv4\sum_{j=k+1}^{\frac{p-1}2}\frac{1}{(2j-1)^2}\\
&=4\left(H_{p-1}^{(2)}-\frac14H_{\frac{p-1}2}^{(2)}-H_{2k}^{(2)}+\frac14H_k^{(2)}\right)\equiv H_{k}^{(2)}-4H_{2k}^{(2)}\pmod{p}.
\end{align*}
\end{proof}
Now we evaluate $\sum_{k=0}^{p-1}\binom{2k}k^2/(-16)^k$ modulo $p^3$.
$$\sum_{k=0}^{p-1}\frac{\binom{2k}k^2}{(-16)^k}=\sum_{k=0}^{(p-1)/2}\frac{\binom{2k}k^2}{(-16)^k}+\sum_{k=(p+1)/2}^{p-1}\frac{\binom{2k}k^2}{(-16)^k}.$$
In view of Lemma \ref{Lem2.1}, we have
\begin{align*}
&\sum_{k=0}^{(p-1)/2}\frac{\binom{2k}k^2}{(-16)^k}\\
&\equiv\sum_{k=0}^{(p-1)/2}\frac{(-1)^k\binom{(p-1)/2}k^2}{\left(1-p\sum_{j=1}^k\frac1{2j-1}+\frac{p^2}2\left(\sum_{j=1}^k\frac1{2j-1}\right)^2-\frac{p^2}2\sum_{j=1}^k\frac1{(2j-1)^2}\right)^2}\\
&\equiv S_1+2pS_2+2p^2S_3+p^2S_4\pmod{p^3},
\end{align*}
where
\begin{gather*}
S_1=\sum_{k=0}^{(p-1)/2}(-1)^k\binom{(p-1)/2}k^2,\\
S_2=\sum_{k=0}^{(p-1)/2}(-1)^k\binom{(p-1)/2}k^2\sum_{j=1}^k\frac1{2j-1},\\
S_3=\sum_{k=0}^{(p-1)/2}(-1)^k\binom{(p-1)/2}k^2\left(\sum_{j=1}^k\frac1{2j-1}\right)^2,\\
S_4=\sum_{k=0}^{(p-1)/2}(-1)^k\binom{(p-1)/2}k^2\sum_{j=1}^k\frac1{(2j-1)^2}.
\end{gather*}
Since $(p-1)/2\equiv1\pmod2$, so we have
$$\sum_{k=0}^{\frac{p-1}2}(-1)^k\binom{\frac{p-1}2}k^2=\sum_{k=0}^{\frac{p-1}2}(-1)^{\frac{p-1}2-k}\binom{\frac{p-1}2}k^2=-\sum_{k=0}^{\frac{p-1}2}(-1)^k\binom{\frac{p-1}2}k^2.$$
Hence
\begin{align}\label{0}
S_1=\sum_{k=0}^{(p-1)/2}(-1)^k\binom{(p-1)/2}k^2=0.
\end{align}
This, with Lemma \ref{Lem2.2} yields that
\begin{align*}
S_2&=\sum_{k=0}^{(p-1)/2}(-1)^k\binom{(p-1)/2}k^2\left(H_{2k}-\frac12H_k\right)\\
&=\sum_{k=0}^{(p-1)/2}(-1)^{(p-1)/2-k}\binom{(p-1)/2}k^2\left(H_{p-1-2k}-\frac12H_{(p-1)/2-k}\right)\\
&\equiv-\sum_{k=0}^{\frac{p-1}2}(-1)^k\binom{(p-1)/2}k^2\left(\frac12H_{k}+\frac{p}4H_k^{(2)}\right)\pmod{p^2}
\end{align*}
and
\begin{align*}
S_3&=\sum_{k=0}^{(p-1)/2}(-1)^k\binom{(p-1)/2}k^2\left(H_{2k}-\frac12H_k\right)^2\\
&=\sum_{k=0}^{(p-1)/2}(-1)^{(p-1)/2-k}\binom{(p-1)/2}k^2\left(H_{p-1-2k}-\frac12H_{(p-1)/2-k}\right)^2\\
&\equiv-\sum_{k=0}^{\frac{p-1}2}(-1)^k\binom{(p-1)/2}k^2\left(\frac14H_{k}^2-\frac{1}2H_kH_{(p-1)/2}\right)\pmod{p}.
\end{align*}
By Lemma \ref{Lem2.2}, we have
\begin{align*}
S_4&=\sum_{k=0}^{(p-1)/2}(-1)^k\binom{(p-1)/2}k^2\left(H_{2k}^{(2)}-\frac14H_k^{(2)}\right)\\
&=\sum_{k=0}^{(p-1)/2}(-1)^{(p-1)/2-k}\binom{(p-1)/2}k^2\left(H_{p-1-2k}^{(2)}-\frac14H_{(p-1)/2-k}^{(2)}\right)\\
&\equiv\frac14\sum_{k=0}^{(p-1)/2}(-1)^{k}\binom{(p-1)/2}k^2H_k^{(2)}\pmod p.
\end{align*}
Therefore,
\begin{align*}
&\sum_{k=0}^{\frac{p-1}2}\frac{\binom{2k}k^2}{(-16)^k}\equiv-p\sum_{k=0}^{\frac{p-1}2}(-1)^k\binom{\frac{p-1}2}k^2H_{k}-\frac{p^2}2\sum_{k=0}^{\frac{p-1}2}(-1)^k\binom{\frac{p-1}2}k^2H_{k}^2\\
&+p^2\sum_{k=0}^{\frac{p-1}2}(-1)^k\binom{\frac{p-1}2}k^2H_kH_{\frac{p-1}2}-\frac{p^2}4\sum_{k=0}^{\frac{p-1}2}(-1)^{k}\binom{\frac{p-1}2}k^2H_k^{(2)}\\
&=-p\sum_{k=0}^{\frac{p-1}2}(-1)^k\binom{\frac{p-1}2}k^2H_{k}+p^2\sum_{k=0}^{\frac{p-1}2}(-1)^k\binom{\frac{p-1}2}k^2H_kH_{\frac{p-1}2}\\
&-\frac{p^2}4\sum_{k=0}^{\frac{p-1}2}(-1)^{k}\binom{\frac{p-1}2}k^2\left(2H_k^2+H_k^{(2)}\right)\pmod{p^3}.
\end{align*}
By \texttt{Sigma}, we can find and prove the following identities,
\begin{align*}
&\sum_{k=0}^{2n-1}(-1)^k\binom{2n-1}k^2H_k=\frac{(-16)^n}{8n\binom{2n}n},\\
&\sum_{k=0}^{2n-1}(-1)^k\binom{2n-1}k^2\left(2H_k^2+H_k^{(2)}\right)\\
&=\frac1{8n^2}-\frac{(-16)^n}{8n^2\binom{2n}n}\bigg(1+2nH_n-6nH_{2n}+n\sum_{k=1}^n\frac{\binom{2k}k}{k(-16)^k}\\
&-12n\sum_{k=1}^n\frac{\binom{2k}k}{(2k-1)^2(-16)^k}-12n\sum_{k=1}^n\frac{\binom{2k}k}{(2k-1)(-16)^k}\bigg).
\end{align*}
Substituting $n=(p+1)/4$ into the above identities, we have
\begin{align*}
\sum_{k=0}^{\frac{p-1}2}(-1)^k\binom{\frac{p-1}2}k^2H_{k}&=\frac{(-1)^{\frac{p+1}4}}{\binom{\frac{p-1}2}{\frac{p-3}4}}\frac{2^{p-1}}{p+1}\\
&\equiv\frac{(-1)^{\frac{p+1}4}(1-p+pq_p(2))}{\binom{\frac{p-1}2}{\frac{p-3}4}}\pmod{p^2}
\end{align*}
and
\begin{align*}
&\sum_{k=0}^{\frac{p-1}2}(-1)^{k}\binom{\frac{p-1}2}k^2\left(2H_k^2+H_k^{(2)}\right)\\
\equiv&2-\frac{4(-1)^{\frac{p+1}4}}{\binom{\frac{p-1}2}{\frac{p-3}4}}\bigg(\frac32q_p(2)+\frac14\sum_{k=1}^{\frac{p+1}4}\frac{\binom{2k}k}{k(-16)^k}-3\sum_{k=1}^{\frac{p+1}4}\frac{\binom{2k}k}{(2k-1)^2(-16)^k}\\
&-3\sum_{k=1}^{\frac{p+1}4}\frac{\binom{2k}k}{(2k-1)(-16)^k}\bigg)\pmod p.
\end{align*}
Therefore, modulo $p^3$, we have
\begin{align*}
&\sum_{k=0}^{\frac{p-1}2}\frac{\binom{2k}k^2}{(-16)^k}\equiv-\frac{p^2}2-\frac{(-1)^{\frac{p+1}4}}{\binom{\frac{p-1}2}{\frac{p-3}4}}p\left(1-p+\frac32pq_p(2)\right)+\frac{p^2}4\frac{(-1)^{\frac{p+1}4}}{\binom{\frac{p-1}2}{\frac{p-3}4}}\times\\
&\bigg(\sum_{k=1}^{\frac{p+1}4}\frac{\binom{2k}k}{k(-16)^k}-12\sum_{k=1}^{\frac{p+1}4}\frac{\binom{2k}k}{(2k-1)^2(-16)^k}-12\sum_{k=1}^{\frac{p+1}4}\frac{\binom{2k}k}{(2k-1)(-16)^k}\bigg).
\end{align*}
\begin{lemma}\label{Lemszw}{\rm (\cite[Lemma 2.1]{sscm})} For any prime $p>3$ and $1\leq j\leq p-1$, we have
\begin{equation*}
j\binom{2j}j\binom{2(p-j)}{p-j}\equiv2p(-1)^{\lfloor 2j/p\rfloor-1}\pmod{p^2}.
\end{equation*}
\end{lemma}
\noindent By this Lemma, we have
\begin{align*}
\sum_{k=\frac{p+1}2}^{p-1}\frac{\binom{2k}k^2}{(-16)^k}&\equiv\sum_{k=\frac{p+1}2}^{p-1}\frac{4p^2}{k^2\binom{2p-2k}{p-k}^2(-16)^k}=\frac{4p^2}{(-16)^p}\sum_{k=1}^{\frac{p-1}2}\frac{(-16)^k}{(p-k)^2\binom{2k}k^2}\\
&\equiv-\frac{p^2}4\sum_{k=1}^{\frac{p-1}2}\frac{(-16)^k}{k^2\binom{2k}k^2}\equiv-\frac{p^2}4\sum_{k=1}^{\frac{p-1}2}\frac{(-1)^k}{k^2\binom{\frac{p-1}2}k^2}\\
&\equiv-p^2\sum_{k=1}^{\frac{p-1}2}\frac{(-1)^k}{\binom{\frac{p-3}2}{k-1}^2}=p^2\sum_{k=0}^{\frac{p-3}2}\frac{(-1)^k}{\binom{\frac{p-3}2}{k}^2}\pmod{p^3}.
\end{align*}
By \texttt{Sigma}, we can find and prove the following identity,
$$
\sum_{k=0}^{2n}\frac{(-1)^k}{\binom{2n}k^2}=(2n+1)\frac{(-1)^n}{\binom{2n}n}\left(\frac34\sum_{k=1}^n\frac{(-1)^k\binom{2k}k}{k}+\sum_{k=0}^n\frac{(-1)^k\binom{2k}k}{2k+1}\right).
$$
Setting $n=(p-3)/4$ into this identity, we have
\begin{align*}
&\sum_{k=0}^{\frac{p-3}2}\frac{(-1)^k}{\binom{\frac{p-3}2}{k}^2}\equiv-\frac12\frac{(-1)^{\frac{p-3}4}}{\binom{\frac{p-3}2}{\frac{p-3}4}}\left(\frac34\sum_{k=1}^{\frac{p-3}4}\frac{(-1)^k\binom{2k}k}{k}+\sum_{k=0}^{\frac{p-3}4}\frac{(-1)^k\binom{2k}k}{2k+1}\right)\\
&\equiv-\frac{(-1)^{\frac{p+1}4}}{\binom{\frac{p-1}2}{\frac{p-3}4}}\left(\frac34\sum_{k=1}^{\frac{p-3}4}\frac{(-1)^k\binom{2k}k}{k}+\sum_{k=0}^{\frac{p-3}4}\frac{(-1)^k\binom{2k}k}{2k+1}\right)\pmod p.
\end{align*}
Therefore,
\begin{align*}
&\sum_{k=0}^{p-1}\frac{\binom{2k}k^2}{(-16)^k}\equiv-\frac{p^2}2-\frac{(-1)^{\frac{p+1}4}}{\binom{\frac{p-1}2}{\frac{p-3}4}}p\left(1-p+\frac32pq_p(2)\right)+\frac{p^2}4\frac{(-1)^{\frac{p+1}4}}{\binom{\frac{p-1}2}{\frac{p-3}4}}\times\\
&\bigg(\sum_{k=1}^{\frac{p+1}4}\frac{\binom{2k}k}{k(-16)^k}-12\sum_{k=1}^{\frac{p+1}4}\frac{\binom{2k}k}{(2k-1)^2(-16)^k}-12\sum_{k=1}^{\frac{p+1}4}\frac{\binom{2k}k}{(2k-1)(-16)^k}\bigg)\\
&-p^2\frac{(-1)^{\frac{p+1}4}}{\binom{\frac{p-1}2}{\frac{p-3}4}}\left(\frac34\sum_{k=1}^{\frac{p-3}4}\frac{(-1)^k\binom{2k}k}{k}+\sum_{k=0}^{\frac{p-3}4}\frac{(-1)^k\binom{2k}k}{2k+1}\right)\pmod{p^3}.
\end{align*}
In view of (\ref{zhu1}), to prove Theorem \ref{Th1.1}, we just need to prove that
\begin{align}\label{zuihou1}
&-2q_p(2)-2\sum_{k=0}^{\frac{p-3}4}\frac{\binom{\frac{p-1}2}k(-1)^k}{2k+1}\equiv(-1)^{\frac{p+1}4}\binom{\frac{p-1}2}{\frac{p-3}4}-\frac12\bigg(\sum_{k=1}^{\frac{p+1}4}\frac{\binom{2k}k}{k(-16)^k}\notag\\
&-12\sum_{k=1}^{\frac{p+1}4}\frac{\binom{2k}k}{(2k-1)^2(-16)^k}-12\sum_{k=1}^{\frac{p+1}4}\frac{\binom{2k}k}{(2k-1)(-16)^k}\bigg)\notag\\
&+2\left(\frac34\sum_{k=1}^{\frac{p-3}4}\frac{(-1)^k\binom{2k}k}{k}+\sum_{k=0}^{\frac{p-3}4}\frac{(-1)^k\binom{2k}k}{2k+1}\right)\pmod p.
\end{align}
While
\begin{align*}
-\frac12\sum_{k=1}^{\frac{p+1}4}\frac{\binom{2k}k}{k(-16)^k}&\equiv-\frac12\sum_{k=1}^{\frac{p+1}4}\frac{\binom{\frac{p-1}2}k}{k4^k}=-\frac12\sum_{k=\frac{p-3}4}^{\frac{p-3}2}\frac{\binom{\frac{p-1}2}k4^k}{((p-1)/2-k)2^{p-1}}\\
&\equiv\sum_{k=\frac{p-3}4}^{\frac{p-3}2}\frac{\binom{\frac{p-1}2}k4^k}{2k+1}\pmod p.
\end{align*}
Hence
\begin{align*}
&-\frac12\sum_{k=1}^{\frac{p+1}4}\frac{\binom{2k}k}{k(-16)^k}+\sum_{k=0}^{\frac{p-3}4}\frac{(-1)^k\binom{2k}k}{2k+1}\equiv\sum_{k=\frac{p-3}4}^{\frac{p-3}2}\frac{\binom{\frac{p-1}2}k4^k}{2k+1}+\sum_{k=0}^{\frac{p-3}4}\frac{(-1)^k\binom{2k}k}{2k+1}\\
&\equiv\sum_{k=\frac{p-3}4}^{\frac{p-3}2}\frac{\binom{2k}k(-1)^k}{2k+1}+\sum_{k=0}^{\frac{p-3}4}\frac{(-1)^k\binom{2k}k}{2k+1}=\sum_{k=0}^{\frac{p-3}2}\frac{\binom{2k}k(-1)^k}{2k+1}-2\binom{\frac{p-3}2}{\frac{p-3}4}(-1)^{\frac{p-3}4}\\
&\equiv\sum_{k=0}^{\frac{p-3}2}\frac{\binom{2k}k(-1)^k}{2k+1}-\binom{\frac{p-1}2}{\frac{p-3}4}(-1)^{\frac{p+1}4}\pmod p.
\end{align*}
So (\ref{zuihou1}) is equivalent to
\begin{align}\label{zuihou2}
&-2q_p(2)-2\sum_{k=0}^{\frac{p-3}4}\frac{\binom{\frac{p-1}2}k(-1)^k}{2k+1}\notag\\
&\equiv\sum_{k=0}^{\frac{p-3}2}\frac{\binom{2k}k(-1)^k}{2k+1}+6\bigg(\sum_{k=1}^{\frac{p+1}4}\frac{\binom{2k}k}{(2k-1)^2(-16)^k}+\sum_{k=1}^{\frac{p+1}4}\frac{\binom{2k}k}{(2k-1)(-16)^k}\bigg)\notag\\
&+\frac32\sum_{k=1}^{\frac{p-3}4}\frac{(-1)^k\binom{2k}k}{k}+\sum_{k=0}^{\frac{p-3}4}\frac{(-1)^k\binom{2k}k}{2k+1}\pmod p.
\end{align}
It is easy to see that
\begin{align*}
&\sum_{k=1}^{\frac{p+1}4}\frac{\binom{2k}k}{(2k-1)^2(-16)^k}=\sum_{k=0}^{\frac{p-3}4}\frac{\binom{2k+2}{k+1}}{(2k+1)^2(-16)^{k+1}}\\
&=-\frac18\sum_{k=0}^{\frac{p-3}4}\frac{\binom{2k}k}{(2k+1)(k+1)(-16)^k}\\
&=-\frac14\sum_{k=0}^{\frac{p-3}4}\frac{\binom{2k}k}{(2k+1)(-16)^k}+\frac18\sum_{k=0}^{\frac{p-3}4}\frac{\binom{2k}k}{(k+1)(-16)^k}
\end{align*}
and
\begin{align*}
&\sum_{k=1}^{\frac{p+1}4}\frac{\binom{2k}k}{(2k-1)(-16)^k}=\sum_{k=0}^{\frac{p-3}4}\frac{\binom{2k+2}{k+1}}{(2k+1)(-16)^{k+1}}\\
&=-\frac18\sum_{k=0}^{\frac{p-3}4}\frac{\binom{2k}k}{(k+1)(-16)^k}.
\end{align*}
Thus,
$$
\sum_{k=1}^{\frac{p+1}4}\frac{\binom{2k}k}{(2k-1)^2(-16)^k}+\sum_{k=1}^{\frac{p+1}4}\frac{\binom{2k}k}{(2k-1)(-16)^k}=-\frac14\sum_{k=0}^{\frac{p-3}4}\frac{\binom{2k}k}{(2k+1)(-16)^k}.
$$
So (\ref{zuihou2}) is equivalent to
\begin{align}\label{zuihou3}
-2q_p(2)-2&\sum_{k=0}^{\frac{p-3}4}\frac{\binom{\frac{p-1}2}k(-1)^k}{2k+1}\equiv\sum_{k=0}^{\frac{p-3}2}\frac{\binom{2k}k(-1)^k}{2k+1}-\frac32\sum_{k=0}^{\frac{p-3}4}\frac{\binom{2k}k}{(2k+1)(-16)^k}\notag\\
&+\frac32\sum_{k=1}^{\frac{p-3}4}\frac{(-1)^k\binom{2k}k}{k}+\sum_{k=0}^{\frac{p-3}4}\frac{(-1)^k\binom{2k}k}{2k+1}\pmod p.
\end{align}
It is easy to get that
\begin{align*}
&\sum_{k=0}^{\frac{p-3}4}\frac{\binom{2k}k}{(2k+1)(-16)^k}\equiv\sum_{k=0}^{\frac{p-3}4}\frac{\binom{\frac{p-1}2}k}{(2k+1)4^k}=\sum_{k=\frac{p+1}4}^{\frac{p-1}2}\frac{\binom{\frac{p-1}2}k4^k}{(p-2k)2^{p-1}}\\
&\equiv-\frac12\sum_{k=\frac{p+1}4}^{\frac{p-1}2}\frac{\binom{\frac{p-1}2}k4^k}{k}\equiv-\frac12\sum_{k=\frac{p+1}4}^{\frac{p-1}2}\frac{\binom{2k}k(-1)^k}{k}\pmod p.
\end{align*}
By \texttt{Sigma}, we can find and prove the following identity,
$$
\sum_{k=0}^n\frac{\binom{2n+1}k(-1)^k}{k-(2n+1)}=-\frac34\sum_{k=1}^n\frac{\binom{2k}k(-1)^k}{k}-\sum_{k=0}^n\frac{\binom{2k}k(-1)^k}{2k+1}.
$$
Substituting $n=(p-3)/4$ into the above identity, we have
\begin{align*}
2\sum_{k=0}^{\frac{p-3}4}\frac{\binom{\frac{p-1}2}k(-1)^k}{2k+1}&\equiv\sum_{k=0}^{\frac{p-3}4}\frac{\binom{\frac{p-1}2}k(-1)^k}{k-(p-1)/2}\\
&=-\frac34\sum_{k=1}^{\frac{p-3}4}\frac{\binom{2k}k(-1)^k}{k}-\sum_{k=0}^{\frac{p-3}4}\frac{\binom{2k}k(-1)^k}{2k+1}
\pmod p.
\end{align*}
Therefore, (\ref{zuihou3}) is equivalent to
\begin{align}\label{zuihou4}
-2q_p(2)\equiv\sum_{k=0}^{\frac{p-3}2}\frac{\binom{2k}k(-1)^k}{2k+1}+\frac34\sum_{k=1}^{\frac{p-1}2}\frac{(-1)^k\binom{2k}k}{k}\pmod p.
\end{align}
In view of \cite[(1.14)]{ST}, we have
\begin{equation}\label{-1k2kkk}
\sum_{k=1}^{\frac{p-1}2}\frac{(-1)^k\binom{2k}k}{k}\equiv-5\frac{F_{p-\left(\frac{p}5\right)}}{p}\pmod p,
\end{equation}
where $\{F_n\}_{n\in\mathbb{N}}$ is the well-known Fibonacci sequence defined by
$$F_0=0,\ F_1=1,\ \mbox{and}\ F_{n+1}=F_n+F_{n-1}\ \mbox{for}\ \ n=1,2,3,\ldots.$$
Its companion Lucas sequence $\{L_n\}_{n\in\mathbb{N}}$ is given by
$$L_0=2,\ L_1=1,\ \mbox{and}\ L_{n+1}=L_n+L_{n-1}\ \mbox{for}\ \ n=1,2,3,\ldots.$$
It is easy to see that
\begin{align*}
\sum_{k=0}^{\frac{p-3}2}\frac{\binom{2k}k(-1)^k}{2k+1}\equiv\sum_{k=1}^{\frac{p-1}2}\frac{\binom{\frac{p-1}2}k2^{p-1}}{(p-2k)4^k}\equiv-\frac12\sum_{k=1}^{\frac{p-1}2}\frac{\binom{2k}k}{k(-16)^k}\pmod{p}.
\end{align*}
Thus, (\ref{zuihou4}) is equivalent to
\begin{align}\label{zuihou5}
\sum_{k=1}^{\frac{p-1}2}\frac{\binom{2k}k}{k(-16)^k}\equiv-\frac{15}2\frac{F_{p-\left(\frac{p}5\right)}}{p}+4q_p(2)\pmod p.
\end{align}
In view of \cite[(32) and (28)]{RM2018} and the last congruence in \cite[Page 15]{ST}, we have
\begin{align*}
\sum_{k=1}^{\frac{p-1}2}\frac{\binom{2k}k}{k(-16)^k}&\equiv2\frac{-\left(\frac12+\frac{\sqrt{5}}{4}\right)^p-\left(\frac12-\frac{\sqrt{5}}{4}\right)^p+1}{p}\\
&=-2\frac{(2+\sqrt{5})^p+(2-\sqrt{5})^p-4^p}{p4^p}\\
&=-2\frac{L_{3p}-4^p}{p4^p}=-2\frac{(L_p-1)(L_p^2+L_p+4)+4-4^p}{p4^p}\\
&\equiv-3\frac{L_p-1}{p}-2\frac{1-4^{p-1}}{p4^{p-1}}\equiv-\frac{15}2\frac{F_{p-\left(\frac{p}5\right)}}{p}+4q_p(2)\pmod p.
\end{align*}
So we obtain Theorem \ref{Th1.1} for each prime $p>5$ with $p\equiv3\pmod4$, and the case of $p=3$ can be checked directly.
Therefore, the proof of Theorem \ref{Th1.1} is complete.\qed
\section{Proof of Theorem \ref{Th1.2}}
 \setcounter{lemma}{0}
\setcounter{theorem}{0}
\setcounter{corollary}{0}
\setcounter{remark}{0}
\setcounter{equation}{0}
\setcounter{conjecture}{0}
In 2010, J. B. Cosgrave and K. Dilcher \cite{CD} obtained that for any prime $p\equiv1\pmod4$, if $p=x^2+y^2$ with $x\equiv1\pmod4$, then
\begin{align}\label{p-12p-14}
\binom{\frac{p-1}2}{\frac{p-1}4}&\equiv\left(2x-\frac{p}{2x}-\frac{p^3}{8x^3}\right)\times\notag\\
&\left(1+\frac12pq_p(2)+\frac18p^2\left(2E_{p-3}-q^2_p(2)\right)\right)\pmod{p^3}.
\end{align}
And in view of \cite[Theorem 4.11]{KMY}, for any prime $p\equiv1\pmod4$, if $p=x^2+y^2$ with $x\equiv1\pmod4$, we have
\begin{equation}\label{3p-34p-14}
\binom{(3p-3)/4}{(p-1)/4}\equiv(-1)^{\frac{p-1}4}\left(2x-\frac{p}{2x}\right)\left(1-\frac{p}2q_p(2)\right)\pmod{p^2}.
\end{equation}
\begin{lemma}\label{Lem3.1} For any prime $p\equiv1\pmod4$ and $p>5$, if $p=x^2+y^2$ with $x\equiv1\pmod4$, then
\begin{gather*}
\sum_{j=0}^{(p-1)/2}\frac{\binom{(p-1)/2}j(-1)^j}{4j+3}\equiv(-1)^{\frac{p-1}4}\frac{1-2p+\frac{p}2q_p(2)+\frac{p}{4x^2}}{2x}\pmod{p^2},\\
\sum_{j=0}^{(p-1)/2}\frac{\binom{(p-1)/2}j(-1)^jH_j}{4j+3}\equiv-\frac1{2x}(-1)^{\frac{p-1}4}\sum_{k=1}^{\frac{p-1}2}\frac{\binom{-3/4}k(-1)^k}{k}\pmod p.
\end{gather*}
\end{lemma}
\begin{proof} By \texttt{Sigma}, we can find and prove the following identities:
\begin{gather*}
\sum_{k=0}^n\frac{\binom{n}k(-1)^k}{4k+3}=\frac{(-1)^n}{(4n+3)\binom{-3/4}n},\\
\sum_{k=0}^n\frac{\binom{n}k(-1)^kH_k}{4k+3}=-\frac{(-1)^n}{(4n+3)\binom{-3/4}n}\sum_{k=0}^n\frac{\binom{-3/4}k(-1)^k}{k}.
\end{gather*}
Similar as the proof of Lemma \ref{mpt}, we can get that
\begin{align*}
\binom{-3/4}{(p-1)/2}&\equiv\binom{(3p-3)/4}{(p-1)/2}\left(1-\frac{3p}4(H_{(3p-3)/2}-H_{(p-1)/4})\right)\\
&\equiv\binom{(3p-3)/4}{(p-1)/2}\pmod{p^2},
 \end{align*}
so by substituting $n=(p-1)/2$ into the above identities and (\ref{3p-34p-14}), we have
\begin{align*}
\sum_{j=0}^{\frac{p-1}2}\frac{\binom{\frac{p-1}2}j(-1)^j}{4j+3}&=\frac{(-1)^{\frac{p-1}2}}{(1+2p)\binom{-3/4}{\frac{p-1}2}}\equiv\frac{1-2p}{\binom{-3/4}{\frac{p-1}2}}\equiv\frac{(-1)^{\frac{p-1}4}(1-2p)}{\left(2x-\frac{p}{2x}\right)\left(1-\frac{p}2q_p(2)\right)}\\
&\equiv(-1)^{\frac{p-1}4}\frac1{2x}\left(1-2p+\frac{p}2q_p(2)+\frac{p}{4x^2}\right)\pmod{p^2}
\end{align*}
and
\begin{align*}
\sum_{j=0}^{(p-1)/2}\frac{\binom{(p-1)/2}j(-1)^jH_j}{4j+3}\equiv-\frac1{2x}(-1)^{\frac{p-1}4}\sum_{k=1}^{\frac{p-1}2}\frac{\binom{-3/4}k(-1)^k}{k}\pmod p.
\end{align*}
These prove Lemma \ref{Lem3.1}.
\end{proof}
\begin{lemma}\label{Lem3.2} For any prime $p\equiv1\pmod4$ and $p>5$, if $p=x^2+y^2$ with $x\equiv1\pmod4$, then
\begin{align*}
&\frac{p}{2\cdot4^{p-1}}\sum_{j=0}^{\frac{p-1}2}\frac{\binom{2j}j\binom{p-1}{2j}}{(4j+1)4^j}\equiv\frac{p}{4^p}(-1)^{\frac{p-1}4}\binom{\frac{p-1}2}{\frac{p-1}4}\sum_{k=1}^{\frac{p-1}2}\frac{\binom{-1/4}k(-1)^k}{k}+(-1)^{\frac{p-1}4}\\
&\times\bigg(x-\frac{p}{4x}-\frac{p^2}{16x^3}-\frac{3px}2q_p(2)+\frac{3p^2}{8x}q_p(2)+3p^2xq^2_p(2)\bigg)\pmod{p^3}.
\end{align*}
\end{lemma}
\begin{proof} In view of Lemma \ref{Lem2.1}, we have
\begin{align*}
&p\sum_{j=0}^{\frac{p-1}2}\frac{\binom{2j}j\binom{p-1}{2j}}{(4j+1)4^j}\equiv p\sum_{\substack{j=0\\j\neq(p-1)/4}}^{\frac{p-1}2}\frac{\binom{2j}j\binom{p-1}{2j}}{(4j+1)4^j}+\frac{\binom{\frac{p-1}2}{\frac{p-1}4}\binom{p-1}{\frac{p-1}2}}{2^{\frac{p-1}2}}\\
&\equiv p\sum_{\substack{j=0\\j\neq\frac{p-1}4}}^{\frac{p-1}2}\frac{\binom{\frac{p-1}2}j(-1)^j}{4j+1}(1-pH_{2j})\left(1+p\left(H_{2j}-\frac12H_j\right)\right)+\frac{\binom{\frac{p-1}2}{\frac{p-1}4}\binom{p-1}{\frac{p-1}2}}{2^{\frac{p-1}2}}\\
&\equiv p\sum_{j=0}^{\frac{p-1}2}\frac{\binom{\frac{p-1}2}j(-1)^j}{4j+1}\left(1-\frac{p}2H_j\right)+T_1\pmod{p^3},
\end{align*}
where
$$
T_1=\frac{\binom{\frac{p-1}2}{\frac{p-1}4}\binom{p-1}{\frac{p-1}2}}{2^{\frac{p-1}2}}-\binom{\frac{p-1}2}{\frac{p-1}4}(-1)^{\frac{p-1}4}\left(1-\frac{p}2H_{(p-1)/4}\right).
$$
It is well-known the congruence of Morley \cite{Mor}:
\begin{equation}\label{mor}
\binom{p-1}{(p-1)/2}\equiv(-1)^{(p-1)/2}4^{p-1}\pmod{p^3}\ \ \mbox{for}\ \mbox{prime}\ \ p>3.
\end{equation}
And in view of \cite[Lemma 5.3]{mw}, we have
$$
2^{\frac{p-1}2}\equiv\left(\frac2p\right)\left(1+\frac{p}2q_p(2)-\frac{p^2}8q^2_p(2)\right)\pmod{p^3}.
$$
These, with Lemma \ref{Lem2.0}, (\ref{p-12p-14}) and $(-1)^{\frac{p-1}4}=\left(\frac2p\right)$ (since $p\equiv1\pmod4$), yield that
\begin{align*}
T_1&\equiv\binom{\frac{p-1}2}{\frac{p-1}4}\left(2^{\frac{3p-3}2}-(-1)^{\frac{p-1}4}\left(1-\frac{p}2H_{(p-1)/4}\right)\right)\\
&\equiv\binom{\frac{p-1}2}{\frac{p-1}4}(-1)^{\frac{p-1}4}\left(1+\frac{3p}2q_p(2)+\frac{3p^2}8q^2_p(2)-\left(1-\frac{p}2H_{(p-1)/4}\right)\right)\\
&\equiv\binom{\frac{p-1}2}{\frac{p-1}4}(-1)^{\frac{p-1}4}\left(\frac{9p^2}8q^2_p(2)-\frac{p^2}2E_{p-3}\right)\\
&\equiv(-1)^{\frac{p-1}4}\left(\frac{9p^2x}4q^2_p(2)-p^2xE_{p-3}\right)\pmod{p^3}.
\end{align*}
Substituting $n=(p-1)/2$ into (\ref{nk4k+1}), (\ref{nk4k+1hk}) and by Lemma \ref{Lem2.0}, we have
\begin{align*}
&p\sum_{k=0}^{\frac{p-1}2}\frac{\binom{\frac{p-1}2}k(-1)^k}{4k+1}=p\frac{(-1)^{\frac{p-1}2}}{(2p-1)\binom{-1/4}{\frac{p-1}2}}\\
&=\frac{\frac{p-1}2!}{\frac14(\frac14+1)\cdots(\frac{p}4-1)(\frac{p}4+1)\cdots(\frac{p}4+\frac{p-1}4)}\\
&=\frac{\frac{p-1}2!}{\prod_{i=1}^{\frac{p-1}4}\left(\frac{p^2}{16}-i^2\right)}\equiv(-1)^{\frac{p-1}4}\binom{\frac{p-1}2}{\frac{p-1}4}\left(1+\frac{p^2}{16}H_{\frac{p-1}4}^{(2)}\right)\\
&\equiv(-1)^{\frac{p-1}4}\binom{\frac{p-1}2}{\frac{p-1}4}\left(1+\frac{p^2}{4}E_{p-3}\right)\pmod{p^3}
\end{align*}
and
\begin{align*}
-\frac{p^2}2\sum_{k=0}^{\frac{p-1}2}\frac{\binom{\frac{p-1}2}k(-1)^kH_k}{4k+1}\equiv\frac{p}2(-1)^{\frac{p-1}4}\binom{\frac{p-1}2}{\frac{p-1}4}\sum_{k=1}^{\frac{p-1}2}\frac{\binom{-1/4}k(-1)^k}{k}\pmod p.
\end{align*}
Hence,
\begin{align*}
p\sum_{j=0}^{\frac{p-1}2}\frac{\binom{2j}j\binom{p-1}{2j}}{(4j+1)4^j}\equiv&(-1)^{\frac{p-1}4}\binom{\frac{p-1}2}{\frac{p-1}4}\left(1+\frac{p^2}4E_{p-3}\right)\\
&+(-1)^{\frac{p-1}4}\left(\frac{9p^2x}4q^2_p(2)-p^2xE_{p-3}\right)\\
&+\frac{p}{4^p}(-1)^{\frac{p-1}4}\binom{\frac{p-1}2}{\frac{p-1}4}\sum_{k=1}^{\frac{p-1}2}\frac{\binom{-1/4}k(-1)^k}{k}\pmod{p^3}.
\end{align*}
Noting that $4^{p-1}=1+2pq_p(2)+p^2q^2_p(2)$ and in view of (\ref{p-12p-14}), we immediately get the desired result. Now the proof of Lemma \ref{Lem3.2} is finished.
\end{proof}
\noindent Now we evaluate $\sum_{k=0}^{p-1}(k+1)\binom{2k}k^2/8^k$ modulo $p^3$. By (\ref{id0}), we have
\begin{align*}
\sum_{k=0}^{p-1}\frac{(k+1)\binom{2k}k^2}{8^k}&=\sum_{k=0}^{p-1}\frac{(k+1)\binom{2k}k}{8^k}\sum_{j=0}^{\lfloor k/2\rfloor}\binom{2j}j\binom{k}{2j}2^{k-2j}\\
&=\sum_{j=0}^{(p-1)/2}\frac{\binom{2j}j}{4^j}\sum_{k=2j}^{p-1}\frac{(k+1)\binom{2k}k\binom{k}{2j}}{4^k}.
\end{align*}
By \texttt{Sigma}, we can find and prove the following identity:
$$
\sum_{k=2j}^{n-1}\frac{(k+1)\binom{2k}k\binom{k}{2j}}{4^k}=\frac{n\binom{2n-1}{n-1}\binom{n-1}{2j}\left(n(4j+1)+4j+2\right)}{4^{n-1}(4j+1)(4j+3)}.
$$
Substituting $n=p$ into the above identity and by (\ref{wol}), we have
\begin{align*}
&\sum_{k=0}^{p-1}\frac{(k+1)\binom{2k}k^2}{8^k}=\sum_{j=0}^{\frac{p-1}2}\frac{\binom{2j}j}{4^j}\frac{p\binom{2p-1}{p-1}\binom{p-1}{2j}\left(p(4j+1)+4j+2\right)}{4^{p-1}(4j+1)(4j+3)}\\
&\equiv \frac{p}{4^{p-1}}\sum_{j=0}^{\frac{p-1}2}\frac{\binom{2j}j\binom{p-1}{2j}}{4^j}\left(\frac{p}{4j+3}+\frac12\frac1{4j+3}+\frac12\frac1{4j+1}\right)\pmod{p^3},
\end{align*}
In view of Lemma \ref{Lem2.1}, we have
\begin{align*}
&\sum_{j=0}^{\frac{p-1}2}\frac{\binom{2j}j(1-pH_{2j})}{(4j+3)4^j}\\
&\equiv\sum_{j=0}^{\frac{p-1}2}\frac{\binom{\frac{p-1}2}j(-1)^j}{4j+3}\left(1+p\left(H_{2j}-\frac12H_j\right)\right)(1-pH_{2j})\\
&\equiv\sum_{j=0}^{\frac{p-1}2}\frac{\binom{\frac{p-1}2}j(-1)^j}{4j+3}\left(1-\frac{p}2H_j\right)\pmod{p^2}.
\end{align*}
This, with Lemma \ref{Lem3.1} yields that
\begin{align*}
&p^2\sum_{k=0}^{p-1}\frac{\binom{2k}k}{(4k+3)4^k}+\frac{p}{2\cdot4^{p-1}}\sum_{k=0}^{\frac{p-1}2}\frac{\binom{2k}k(1-pH_{2k})}{(4k+3)4^k}\\
&\equiv \frac{p^2}{2x}(-1)^{\frac{p-1}4}+\frac{p}{2\cdot4^{p-1}}\sum_{j=0}^{\frac{p-1}2}\frac{\binom{\frac{p-1}2}j(-1)^j}{(4j+3)}\left(1-\frac{p}2H_j\right)\\
&\equiv(-1)^{\frac{p-1}4}\left(\frac{p}{4x}-\frac{3p^2}{8x}q_p(2)+\frac{p^2}{16x^3}+\frac{p^2}{8x}\sum_{k=1}^{\frac{p-1}2}\frac{\binom{-3/4}k(-1)^k}{k}\right)\pmod{p^3}.
\end{align*}
Hence by Lemma \ref{Lem3.2}, we can obtain that
\begin{align*}
\sum_{k=0}^{p-1}\frac{(k+1)\binom{2k}k^2}{8^k}\equiv&(-1)^{\frac{p-1}4}\left(x-\frac{3px}2q_p(2)+3p^2xq^2_p(2)\right)\\
&+(-1)^{\frac{p-1}4}\frac{p^2}{8x}\sum_{k=1}^{\frac{p-1}2}\frac{\binom{-3/4}k(-1)^k}{k}\\
&+\frac{p}{4^p}(-1)^{\frac{p-1}4}\binom{\frac{p-1}2}{\frac{p-1}4}\sum_{k=1}^{\frac{p-1}2}\frac{\binom{-1/4}k(-1)^k}{k}\pmod {p^3}.
\end{align*}
In view of \cite[Lemma 3.1]{sijnt} and Lemma \ref{Lem2.0}, we have
\begin{align*}
\sum_{k=1}^{\frac{p-1}2}\frac{\binom{-\frac14}k(-1)^k}{k}&\equiv\frac{p}4H_{\frac{p-1}2}^{(2)}-H_{\frac{p-1}2}-\frac{p}4H_{\frac{p-1}4}^{(2)}-\frac{p}4\sum_{k=1}^{\frac{p-1}4}\frac{(-1)^k}{k^2\binom{\frac{p-1}2}{k}}\\
&\equiv3q_p(2)-\frac{3p}2q^2_p(2)-\frac{p}4\sum_{k=1}^{\frac{p-1}4}\frac{(-1)^k}{k^2\binom{\frac{p-1}2}{k}}\pmod{p^2}.
\end{align*}
This, with (\ref{p-12p-14}) and $4^{p-1}=1+2pq_p(2)+p^2q^2_p(2)$, yields that
\begin{align*}
&\frac{p}{4^p}(-1)^{\frac{p-1}4}\binom{\frac{p-1}2}{\frac{p-1}4}\sum_{k=1}^{\frac{p-1}2}\frac{\binom{-1/4}k(-1)^k}{k}\equiv(-1)^{\frac{p-1}4}\times\\
&\left(\frac{3px}2q_p(2)-\frac{3p^2}{8x}q_p(2)-3p^2xq^2_p(2)-\frac{p^2x}8\sum_{k=1}^{\frac{p-1}4}\frac{(-1)^k}{k^2\binom{\frac{p-1}2}{k}}\right)\pmod{p^3}.
\end{align*}
In view of \cite[(1.45), (1.134)]{G}, we have
\begin{align*}
&\sum_{k=1}^{\frac{p-1}2}\frac{\binom{-3/4}k(-1)^k}{k}\equiv\sum_{k=1}^{\frac{p-1}2}\frac{\binom{\frac{3p-3}4}k(-1)^k}{k}\\
&=\sum_{k=1}^{\frac{3p-3}4}\frac{\binom{\frac{3p-3}4}k(-1)^k}{k}-\sum_{k=\frac{p+1}2}^{\frac{3p-3}4}\frac{\binom{\frac{3p-3}4}k(-1)^k}{k}\\
&=-H_{\frac{3p-3}4}+\sum_{k=\frac{p+1}2}^{\frac{3p-3}4}\binom{k-1}{\frac{p-1}2}\frac{1}{k}\\
&\equiv3q_p(2)+2\sum_{k=0}^{\frac{p-5}4}\frac{\binom{\frac{p-1}2}k(-1)^k}{2k+1}\pmod p.
\end{align*}
Therefore,
\begin{align}\label{k+18k}
&\sum_{k=0}^{p-1}\frac{(k+1)\binom{2k}k^2}{8^k}\equiv(-1)^{\frac{p-1}4}\times\notag\\
&\left(x-\frac{p^2x}8\sum_{k=1}^{\frac{p-1}4}\frac{(-1)^k}{k^2\binom{\frac{p-1}2}{k}}+\frac{p^2}{4x}\sum_{k=0}^{\frac{p-5}4}\frac{\binom{\frac{p-1}2}k(-1)^k}{2k+1}\right)\pmod{p^3}.
\end{align}
Now we evaluate $\sum_{k=0}^{(p-1)/2}(2k+1)\binom{2k}k^2/(-16)^k$ modulo $p^3$. In view of Lemma \ref{Lem2.1}, we have
\begin{align*}
&\sum_{k=0}^{\frac{p-1}2}\frac{(2k+1)\binom{2k}k^2}{(-16)^k}\\
&\equiv\sum_{k=0}^{\frac{p-1}2}\frac{(2k+1)(-1)^k\binom{(p-1)/2}k^2}{\left(1-p\sum_{j=1}^k\frac1{2j-1}+\frac{p^2}2\left(\sum_{j=1}^k\frac1{2j-1}\right)^2-\frac{p^2}2\sum_{j=1}^k\frac1{(2j-1)^2}\right)^2}\\
&\equiv A_1+2pA_2+2p^2A_3+p^2A_4\pmod{p^3},
\end{align*}
where
\begin{gather*}
A_1=\sum_{k=0}^{(p-1)/2}(2k+1)(-1)^k\binom{(p-1)/2}k^2,\\
A_2=\sum_{k=0}^{(p-1)/2}(2k+1)(-1)^k\binom{(p-1)/2}k^2\left(H_{2k}-\frac12H_k\right),\\
A_3=\sum_{k=0}^{(p-1)/2}(2k+1)(-1)^k\binom{(p-1)/2}k^2\left(H_{2k}-\frac12H_k\right)^2,\\
A_4=\sum_{k=0}^{(p-1)/2}(2k+1)(-1)^k\binom{(p-1)/2}k^2\left(H_{2k}^{(2)}-\frac14H_k^{(2)}\right).
\end{gather*}
By Lemma \ref{Lem2.2}, we have
\begin{align*}
A_2&=\sum_{k=0}^{\frac{p-1}2}(2k+1)(-1)^k\binom{\frac{p-1}2}k^2\left(H_{2k}-\frac12H_k\right)\\
&=\sum_{k=0}^{\frac{p-1}2}(p-2k)(-1)^{\frac{p-1}2-k}\binom{\frac{p-1}2}k^2\left(H_{p-1-2k}-\frac12H_{\frac{p-1}2-k}\right)\\
&\equiv\sum_{k=0}^{\frac{p-1}2}(p-2k)(-1)^k\binom{\frac{p-1}2}k^2\left(\frac12H_{k}-\frac12H_{\frac{p-1}2}+\frac{p}4H_k^{(2)}\right)\pmod{p^2}
\end{align*}
and
\begin{align*}
A_3&=\sum_{k=0}^{\frac{p-1}2}(2k+1)(-1)^k\binom{\frac{p-1}2}k^2\left(H_{2k}-\frac12H_k\right)^2\\
&=\sum_{k=0}^{\frac{p-1}2}(p-2k)(-1)^{\frac{p-1}2-k}\binom{\frac{p-1}2}k^2\left(H_{p-1-2k}-\frac12H_{\frac{p-1}2-k}\right)^2\\
&\equiv-\frac12\sum_{k=0}^{\frac{p-1}2}k(-1)^k\binom{\frac{p-1}2}k^2\left(H_{k}^2-2H_kH_{\frac{p-1}2}+H_{\frac{p-1}2}^2\right)\pmod{p}.
\end{align*}
By Lemma \ref{Lem2.2}, we have
\begin{align*}
A_4&=\sum_{k=0}^{\frac{p-1}2}(2k+1)(-1)^k\binom{\frac{p-1}2}k^2\left(H_{2k}^{(2)}-\frac14H_k^{(2)}\right)\\
&=\sum_{k=0}^{\frac{p-1}2}(p-2k)(-1)^{\frac{p-1}2-k}\binom{\frac{p-1}2}k^2\left(H_{p-1-2k}^{(2)}-\frac14H_{\frac{p-1}2-k}^{(2)}\right)\\
&\equiv\frac12\sum_{k=0}^{\frac{p-1}2}k(-1)^{k}\binom{\frac{p-1}2}k^2H_k^{(2)}\pmod p.
\end{align*}
Therefore,
\begin{align*}
&\sum_{k=0}^{\frac{p-1}2}\frac{(2k+1)\binom{2k}k^2}{(-16)^k}\\
\equiv&\sum_{k=0}^{\frac{p-1}2}(2k+1)(-1)^k\binom{\frac{p-1}2}k^2+p^2\sum_{k=0}^{\frac{p-1}2}(-1)^k\binom{\frac{p-1}2}k^2(H_{k}-H_{\frac{p-1}2})\\
&-2p\sum_{k=0}^{\frac{p-1}2}k(-1)^k\binom{\frac{p-1}2}k^2(H_{k}-H_{\frac{p-1}2})\\
&+p^2\sum_{k=0}^{\frac{p-1}2}k(-1)^k\binom{\frac{p-1}2}k^2(2H_kH_{\frac{p-1}2}-H_{\frac{p-1}2}^2)\\
&-\frac{p^2}2\sum_{k=0}^{\frac{p-1}2}k(-1)^{k}\binom{\frac{p-1}2}k^2(2H_k^2+H_k^{(2)})\pmod{p^3}.
\end{align*}
By \texttt{Sigma}, we can find and prove the following identities,
\begin{align*}
&\sum_{k=0}^{2n}(2k+1)(-1)^k\binom{2n}k^2=(-1)^n(2n+1)\binom{2n}n,\\
&\sum_{k=0}^{2n}(-1)^k\binom{2n}k^2=(-1)^n\binom{2n}n,\\
&\sum_{k=0}^{2n}(-1)^k\binom{2n}k^2H_k=\frac12(-1)^n\binom{2n}n(H_{2n}+H_n),\\
&\sum_{k=0}^{2n}(-1)^kk\binom{2n}k^2H_k\\
&=-\frac14\frac{(-16)^n}{\binom{2n}n}+\frac{(-1)^n}4\binom{2n}n+\frac{(-1)^nn}2\binom{2n}n(H_{2n}+H_n).
\end{align*}
Substituting $n=(p-1)/4$ into the above identities, then by Lemma \ref{Lem2.0} and (\ref{p-12p-14}), we have
\begin{align*}
&\sum_{k=0}^{\frac{p-1}2}(2k+1)(-1)^k\binom{\frac{p-1}2}k^2-2p\sum_{k=0}^{\frac{p-1}2}k(-1)^k\binom{\frac{p-1}2}k^2(H_{k}-H_{\frac{p-1}2})\\
\equiv&(-1)^{\frac{p-1}4}\left(x-\frac{p^2x}8q^2_p(2)+\frac{p^2x}2q_p(2)+\frac{p^2}{8x}q_p(2)-\frac{p^2x}4E_{p-3}\right)\pmod{p^3},
\end{align*}
\begin{align*}
p^2\sum_{k=0}^{\frac{p-1}2}(-1)^k\binom{\frac{p-1}2}k^2(H_{k}-H_{\frac{p-1}2})&=\frac{p^2}2(-1)^{\frac{p-1}4}\binom{\frac{p-1}2}{\frac{p-1}4}(H_{\frac{p-1}4}-H_{\frac{p-1}2})\\
&\equiv-p^2x(-1)^{\frac{p-1}4}q_p(2)\pmod{p^3}
\end{align*}
and
\begin{align*}
&p^2\sum_{k=0}^{\frac{p-1}2}k(-1)^k\binom{\frac{p-1}2}k^2(2H_kH_{\frac{p-1}2}-H_{\frac{p-1}2}^2)\\
&\equiv-4p^2q_p(2)\sum_{k=0}^{\frac{p-1}2}k(-1)^k\binom{\frac{p-1}2}k^2H_k-4p^2q^2_p(2)\sum_{k=0}^{\frac{p-1}2}k(-1)^k\binom{\frac{p-1}2}k^2\\
&\equiv(-1)^{\frac{p-1}4}\left(\frac{p^2}{2x}q_p(2)-2p^2xq_p(2)-3p^2xq_p(2)^2\right)\pmod{p^3}.
\end{align*}
Therefore,
\begin{align}\label{2k+1-16k}
&\sum_{k=0}^{\frac{p-1}2}\frac{(2k+1)\binom{2k}k^2}{(-16)^k}\notag\\
&\equiv(-1)^{\frac{p-1}4}\left(x+\frac{5p^2}{8x}q_p(2)-\frac{5p^2x}2q_p(2)-\frac{25p^2x}8q_p(2)^2-\frac{p^2x}4E_{p-3}\right)\notag\\
&\ \ \ \ \ \ \ -\frac{p^2}2\sum_{k=0}^{\frac{p-1}2}k(-1)^{k}\binom{\frac{p-1}2}k^2(2H_k^2+H_k^{(2)})\pmod{p^3}.
\end{align}
By \texttt{Sigma}, we could find and prove the following complicated identity:
\begin{align*}
&\sum_{k=0}^{2n}k(-1)^{k}\binom{2n}k^2(2H_k^2+H_k^{(2)})=\frac{2n-3}4+\frac34\frac{(-16)^n}{\binom{2n}n}-\frac{n}2(-1)^n\binom{2n}n\\
&-\frac32\frac{(-16)^n}{\binom{2n}n}\left(H_{2n}-\frac13H_n\right)+\frac{n(-1)^n}2\binom{2n}n\left(H_{2n}^{(2)}+\frac12H_n^{(2)}\right)\\
&+\frac{(-1)^n}2\binom{2n}n(H_{2n}+H_n)+\frac{n(-1)^n}2\binom{2n}n\left(H_{2n}+H_n\right)^2\\
&-\frac{5n(-1)^n}2\binom{2n}n\sum_{k=1}^n\frac{(-1)^k}{\binom{2k}k}+\frac{3n(-1)^n}4\binom{2n}n\sum_{k=1}^n\frac{(-1)^k}{k^2\binom{2k}k}\\
&+2n(-1)^n\binom{2n}n\sum_{k=1}^n\frac{(-1)^k}{(2k-1)\binom{2k}k}+\frac{15}4\frac{(-16)^n}{\binom{2n}n}\sum_{k=1}^n\frac{\binom{2k}k}{(-16)^k}\\
&+\frac14\frac{(-16)^n}{\binom{2n}n}\sum_{k=1}^n\frac{\binom{2k}k}{k(-16)^k}-3\frac{(-16)^n}{\binom{2n}n}\sum_{k=1}^n\frac{\binom{2k}k}{(2k-1)^2(-16)^k}.
\end{align*}
Setting $n=(p-1)/4$ into this identity, then by Lemma \ref{Lem2.0} and (\ref{p-12p-14}), we have
\begin{align*}
&\sum_{k=0}^{\frac{p-1}2}k(-1)^{k}\binom{\frac{p-1}2}k^2(2H_k^2+H_k^{(2)})\equiv\\
-&\frac78+(-1)^{\frac{p-1}4}\left(\frac{3}{8x}+\frac{x}4-\frac{25x}4q^2_p(2)+\frac{3}{4x}q_p(2)-\frac{x}2E_{p-3}-5xq_p(2)\right)\\
+&(-1)^{\frac{p-1}4}\bigg(\frac{5x}4\sum_{k=1}^{\frac{p-1}4}\frac{(-1)^k}{\binom{2k}k}-\frac{3x}8\sum_{k=1}^{\frac{p-1}4}\frac{(-1)^k}{k^2\binom{2k}k}-x\sum_{k=1}^{\frac{p-1}4}\frac{(-1)^k}{(2k-1)\binom{2k}k}\\
+&\frac{15}{8x}\sum_{k=1}^{\frac{p-1}4}\frac{\binom{2k}k}{(-16)^k}+\frac1{8x}\sum_{k=1}^{\frac{p-1}4}\frac{\binom{2k}k}{k(-16)^k}-\frac3{2x}\sum_{k=1}^{\frac{p-1}4}\frac{\binom{2k}k}{(2k-1)^2(-16)^k}\bigg)\pmod p.
\end{align*}
Substituting this into (\ref{2k+1-16k}), then combining (\ref{k+18k}), we can see that to prove Theorem \ref{Th1.2}, we just need to prove
\begin{align}\label{zhu21}
&-\frac3{16x}-\frac{x}8+\frac{q_p(2)}{4x}+\frac7{16}(-1)^{\frac{p-1}4}-\frac{x}8\sum_{k=1}^{\frac{p-1}4}\frac{(-1)^k}{k^2\binom{\frac{p-1}2}k}+\frac1{4x}\sum_{k=0}^{\frac{p-5}4}\frac{\binom{\frac{p-1}2}k(-1)^k}{2k+1}\notag\\
&-\frac{5x}8\sum_{k=1}^{\frac{p-1}4}\frac{(-1)^k}{\binom{2k}k}+\frac{3x}{16}\sum_{k=1}^{\frac{p-1}4}\frac{(-1)^k}{k^2\binom{2k}k}+\frac{x}2\sum_{k=1}^{\frac{p-1}4}\frac{(-1)^k}{(2k-1)\binom{2k}k}-\frac{15}{16x}\sum_{k=1}^{\frac{p-1}4}\frac{\binom{2k}k}{(-16)^k}\notag\\
&-\frac1{16x}\sum_{k=1}^{\frac{p-1}4}\frac{\binom{2k}k}{k(-16)^k}+\frac3{4x}\sum_{k=1}^{\frac{p-1}4}\frac{\binom{2k}k}{(2k-1)^2(-16)^k}\equiv0\pmod{p}.
\end{align}
By \texttt{Sigma}, we could find and prove the following identity:
\begin{align*}
\sum_{k=1}^n\frac{(-1)^k}{k^2\binom{2n}k}=H_{2n}^{(2)}+\frac32\sum_{k=1}^{n}\frac{(-1)^k}{k^2\binom{2k}k}-2\sum_{k=1}^{n}\frac{(-1)^k}{k\binom{2k}k}+4\sum_{k=1}^{n}\frac{(-1)^k}{(2k-1)\binom{2k}k}.
\end{align*}
Setting $n=(p-1)/4$ into this identity, then combining (\ref{zhu21}), we just need to prove
\begin{align}\label{zhu22}
&-\frac3{16x}-\frac{x}8+\frac{q_p(2)}{4x}+\frac7{16}(-1)^{\frac{p-1}4}+\frac{x}4\sum_{k=1}^{\frac{p-1}4}\frac{(-1)^k}{k\binom{2k}k}+\frac1{4x}\sum_{k=0}^{\frac{p-5}4}\frac{\binom{\frac{p-1}2}k(-1)^k}{2k+1}\notag\\
&-\frac{5x}8\sum_{k=1}^{\frac{p-1}4}\frac{(-1)^k}{\binom{2k}k}-\frac{15}{16x}\sum_{k=1}^{\frac{p-1}4}\frac{\binom{2k}k}{(-16)^k}-\frac1{16x}\sum_{k=1}^{\frac{p-1}4}\frac{\binom{2k}k}{k(-16)^k}\notag\\
&+\frac3{4x}\sum_{k=1}^{\frac{p-1}4}\frac{\binom{2k}k}{(2k-1)^2(-16)^k}\equiv0\pmod{p}.
\end{align}
By \texttt{Sigma}, we could find and prove that
$$
\sum_{k=1}^n\frac{\binom{2n}k(-1)^k}{k}=-H_{2n}+\frac34\sum_{k=1}^n\frac{\binom{2k}k(-1)^k}k-\frac14\sum_{k=1}^n\frac{2k\binom{2k}k(-1)^k}{(2k-1)^2}.
$$
Substituting $n=(p-1)/4$ into this identity, we have
\begin{align*}
&\sum_{k=1}^{\frac{p-1}4}\frac{\binom{\frac{p-1}2}k(-1)^k}{k}=-H_{\frac{p-1}2}+\frac34\sum_{k=1}^{\frac{p-1}4}\frac{\binom{2k}k(-1)^k}k+\sum_{k=0}^{\frac{p-5}4}\frac{\binom{2k}k(-1)^k}{2k+1}\\
&\equiv-H_{\frac{p-1}2}+\frac34\sum_{k=1}^{\frac{p-1}4}\frac{\binom{2k}k(-1)^k}k-\frac12\sum_{k=\frac{p+3}4}^{\frac{p-1}2}\frac{\binom{2k}k}{k(-16)^k}\pmod p.
\end{align*}
Hence,
\begin{align*}
\sum_{k=0}^{\frac{p-5}4}\frac{\binom{\frac{p-1}2}k(-1)^k}{2k+1}&=\sum_{k=\frac{p+3}4}^{\frac{p-1}2}\frac{\binom{\frac{p-1}2}k(-1)^k}{p-2k}\\
&\equiv-\frac12\sum_{k=1}^{\frac{p-1}2}\frac{\binom{\frac{p-1}2}k(-1)^k}k+\frac12\sum_{k=1}^{\frac{p-1}4}\frac{\binom{\frac{p-1}2}k(-1)^k}{k}\\
&\equiv\frac38\sum_{k=1}^{\frac{p-1}4}\frac{\binom{2k}k(-1)^k}k-\frac14\sum_{k=\frac{p+3}4}^{\frac{p-1}2}\frac{\binom{2k}k}{k(-16)^k}\pmod p.
\end{align*}
Thus, (\ref{zhu22}) is equivalent to
\begin{align}\label{zhu23}
&-\frac3{16x}-\frac{x}8+\frac{q_p(2)}{4x}+\frac7{16}(-1)^{\frac{p-1}4}+\frac{x}4\sum_{k=1}^{\frac{p-1}4}\frac{(-1)^k}{k\binom{2k}k}+\frac3{32x}\sum_{k=1}^{\frac{p-1}4}\frac{\binom{2k}k(-1)^k}k\notag\\
&-\frac1{16x}\sum_{k=1}^{\frac{p-1}2}\frac{\binom{2k}k}{k(-16)^k}-\frac{5x}8\sum_{k=1}^{\frac{p-1}4}\frac{(-1)^k}{\binom{2k}k}-\frac{15}{16x}\sum_{k=1}^{\frac{p-1}4}\frac{\binom{2k}k}{(-16)^k}\notag\\
&+\frac3{4x}\sum_{k=1}^{\frac{p-1}4}\frac{\binom{2k}k}{(2k-1)^2(-16)^k}\equiv0\pmod{p}.
\end{align}
It is easy to check that
\begin{align*}
&\sum_{k=1}^{\frac{p-1}4}\frac{\binom{2k}k}{(2k-1)^2(-16)^k}=\sum_{k=0}^{\frac{p-5}4}\frac{2\binom{2k}k}{(2k+1)(k+1)(-16)^{k+1}}\\
&=-\frac18\sum_{k=0}^{\frac{p-5}4}\frac{\binom{2k}k}{(-16)^k}\left(\frac2{2k+1}-\frac1{k+1}\right)\\
&\equiv\frac18\sum_{k=\frac{p+3}4}^{\frac{p-1}2}\frac{\binom{2k}k(-1)^k}{k}+\frac18\sum_{k=0}^{\frac{p-5}4}\frac{\binom{2k}k}{(k+1)(-16)^k}\pmod p.
\end{align*}
This, with (\ref{-1k2kkk}) and (\ref{zuihou5}) yields that (\ref{zhu23}) is equivalent to
\begin{align}\label{zhu24}
&-\frac3{16x}-\frac{x}8+\frac7{16}(-1)^{\frac{p-1}4}+\frac{x}4\sum_{k=1}^{\frac{p-1}4}\frac{(-1)^k}{k\binom{2k}k}-\frac{5x}8\sum_{k=1}^{\frac{p-1}4}\frac{(-1)^k}{\binom{2k}k}\notag\\
&-\frac{15}{16x}\sum_{k=1}^{\frac{p-1}4}\frac{\binom{2k}k}{(-16)^k}+\frac3{32x}\sum_{k=0}^{\frac{p-5}4}\frac{\binom{2k}k}{(k+1)(-16)^k}\equiv0\pmod{p}.
\end{align}
Setting $t=-1, n=(p-1)/4$ into \cite[(25)]{RM} and by (\ref{p-12p-14}), we can obtain that
$$
\frac{x}4\sum_{k=1}^{\frac{p-1}4}\frac{(-1)^k}{k\binom{2k}k}-\frac{5x}8\sum_{k=1}^{\frac{p-1}4}\frac{(-1)^k}{\binom{2k}k}\equiv\frac{x}8-\frac1{16}(-1)^{\frac{p-1}4}\pmod p.
$$
Thus, to prove (\ref{zhu24}), we just need to prove
\begin{align}\label{zhu25}
&-\frac3{16x}+\frac3{8}(-1)^{\frac{p-1}4}-\frac{15}{16x}\sum_{k=1}^{\frac{p-1}4}\frac{\binom{2k}k}{(-16)^k}+\frac3{32x}\sum_{k=0}^{\frac{p-5}4}\frac{\binom{2k}k}{(k+1)(-16)^k}\notag\\
&\equiv0\pmod{p}.
\end{align}
By \texttt{Sigma}, we could find and prove the following identity:
\begin{align*}
\sum_{k=1}^n\frac{\binom{2n}k\left(\frac1{k+1}-10\right)}{4^k}=&-\frac{2n+1}{n+1}\frac{\binom{2n}n}{4^n}+\frac{(18n+5)2^{4n}-(16n+4)5^{2n}}{(2n+1)2^{4n}}\\
&-\frac{3(4n+1)}{2(2n+1)}\frac{5^{2n}}{2^{4n}}\sum_{k=1}^n\binom{2k}k\left(\frac4{25}\right)^k.
\end{align*}
Substituting $n=(p-1)/4$ into this identity and by (\ref{p-12p-14}), we have
\begin{align*}
\sum_{k=1}^{\frac{p-1}4}\frac{\binom{2k}k}{(k+1)(-16)^k}-10\sum_{k=1}^{\frac{p-1}4}\frac{\binom{2k}k}{(-16)^k}&\equiv\sum_{k=1}^{\frac{p-1}4}\frac{\binom{\frac{p-1}2}k\left(\frac1{k+1}-10\right)}{4^k}\\
&\equiv1-\frac{4x}3(-1)^{\frac{p-1}4}\pmod p.
\end{align*}
This, with the fact that
\begin{align*}
\sum_{k=0}^{\frac{p-5}4}\frac{\binom{2k}k}{(k+1)(-16)^k}&=1-\frac{\binom{\frac{p-1}2}{\frac{p-1}4}}{\frac{p+3}4(-16)^{\frac{p-1}4}}+\sum_{k=1}^{\frac{p-1}4}\frac{\binom{2k}k}{(k+1)(-16)^k}\\
&\equiv1-\frac{8x}3(-1)^{\frac{p-1}4}+\sum_{k=1}^{\frac{p-1}4}\frac{\binom{2k}k}{(k+1)(-16)^k}\pmod p
\end{align*}
yields that (\ref{zhu25}) is correct.

\noindent So we obtain Theorem \ref{Th1.2} for each prime $p>5$ with $p\equiv1\pmod4$ and $p=x^2+y^2$ ($4|(x-1)\ \&\ 2|y$), the case of $p=5$ can be checked directly.

\noindent Therefore, the proof of Theorem \ref{Th1.2} is complete.\qed

\Ack. The author would like to thank Prof. R. Tauraso for giving the reference \cite{RM2018} to prove (\ref{zuihou5}). The author was supported by the National Natural Science Foundation of China (grant no. 12001288).

     \end{document}